# A basic class of symmetric orthogonal polynomials using the extended Sturm-Liouville theorem for symmetric functions

**Mohammad Masjed-Jamei**

*Department of Mathematics, K.N.Toosi University of Technology, P.O.Box 16315-1618, Tehran, Iran,*
*E-mail: mmjamei@kntu.ac.ir , mmjamei@yahoo.com*

**Abstract.** In this research, by applying the extended Sturm-Liouville theorem for symmetric functions, a basic class of symmetric orthogonal polynomials (BCSOP) with four free parameters is introduced and all its standard properties, such as a generic second order differential equation along with its explicit polynomial solution, a generic orthogonality relation, a generic three term recurrence relation and so on, are presented. Then, it is shown that four main sequences of symmetric orthogonal polynomials can essentially be extracted from the introduced class. They are respectively the generalized ultraspherical polynomials, generalized Hermite polynomials and two other sequences of symmetric polynomials, which are finitely orthogonal on $(-\infty, \infty)$ and can be expressed in terms of the mentioned class directly. In this way, two half-trigonometric sequences of orthogonal polynomials, as special sub-cases of BCSOP, are also introduced.

**Keywords.** extended Sturm-Liouville theorem for symmetric functions, orthogonal polynomials, Favard's theorem, Pearson distributions family, dual symmetric distributions family, generalized ultraspherical polynomials, fifth and sixth kind of Chebyshev polynomials, generalized Hermite polynomials, two kinds of finite classical symmetric orthogonal polynomials.

**MSC (2000):** 05E05, 33C47, 05E35, 34B24

## 1. Introduction
As it is known, to prove the orthogonality property of a sequence of functions, which satisfies the differential equation

$$p(x) y_n''(x) + q(x) y_n'(x) + (\lambda_n r(x) + s(x)) y_n(x) = 0 , \tag{1}$$

one of the current methods is to apply the Sturm-Liouville theorem ( see e.g. [3] , [19]). According to this theorem, a necessary condition to establish the orthogonality property is that the coefficients $p(x)$, $q(x)$, $r(x)$ and $s(x)$ in (1) must be independent of $n$.

But what will happen if one of these mentioned coefficients depends on *n*? Can we again use the Sturm-Liouville theorem for this case? We have responded to these questions in [16] in detail. However, since the main theorem of [16] is required in this article, we restate it here.

**Theorem 1.** *Let $\Phi_n(x) = (-1)^n \Phi_n(-x)$ be a sequence of symmetric functions that satisfies the differential equation*

$$A(x)\Phi_n''(x) + B(x)\Phi_n'(x) + \left(\lambda_n C(x) + D(x) + (1-(-1)^n)E(x)/2\right)\Phi_n(x) = 0, \qquad (2)$$

*where $A(x), B(x), C(x), D(x)$ and $E(x)$ are independent functions and $\{\lambda_n\}$ is a sequence of constants. If $A(x), C(x), D(x)$ and $E(x)$ are even functions and $B(x)$ is an odd function then*

$$\int_{-\alpha}^{\alpha} W^*(x)\Phi_n(x)\Phi_m(x)\,dx = \left(\int_{-\alpha}^{\alpha} W^*(x)\Phi_n^2(x)\,dx\right)\delta_{n,m}, \qquad (3)$$

*where $\delta_{n,m} = \begin{cases} 0 & \text{if } n \neq m \\ 1 & \text{if } n = m \end{cases}$ and $W^*(x)$ denotes the corresponding weight function as*

$$W^*(x) = C(x)\exp(\int \frac{B(x)-A'(x)}{A(x)}\,dx) = \frac{C(x)}{A(x)}\exp(\int \frac{B(x)}{A(x)}\,dx). \qquad (4)$$

*Of course, $W^*(x)$ defined in (4) must be a positive and even function on $[-\alpha, \alpha]$ and $x = \alpha$ must be a root of the function*

$$A(x)K(x) = A(x)\exp(\int \frac{B(x)-A'(x)}{A(x)}\,dx) = \exp(\int \frac{B(x)}{A(x)}\,dx), \qquad (5)$$

*i.e. $A(\alpha)K(\alpha) = 0$. Note that since*

$$K(x) = \frac{W^*(x)}{C(x)}, \qquad (6)$$

*this function is even ($K(-x) = K(x)$) and therefore $A(-\alpha)K(-\alpha) = 0$ automatically.*

In [16] we have called $\{\Phi_n(x)\}$, "A generalized symmetric sequence of classical orthogonal functions", because by using it, one can generalize several well-known *symmetric* orthogonal functions such as trigonometric orthogonal sequences related to Fourier analysis [1,19], the associated Legendre functions $P_n^m(x)$ [3] having the differential equation

$$(1-x^2)\Phi_n''(x) - 2x\Phi_n'(x) + (n(n+1) - \frac{m^2}{1-x^2})\Phi_n(x) = 0, \qquad (7)$$

and the Bessel functions $J_n(x)$ [3] of symmetric type. For instance, if

$$A(x) = 1 - x^2 = A(-x), \qquad\qquad B(x) = -2x = -B(-x),$$

$$C(x) = 1 = C(-x), \qquad\qquad D(x) = -\frac{m^2}{1-x^2} = D(-x), \qquad (8)$$

$$\lambda_n = n(n+1), \qquad\qquad E(x) = E(-x) \ ; \ Arbitrary,$$

are selected in the main equation (2), the corresponding solution generates the generalized associated Legendre functions and *preserves* the orthogonality property on the same interval [-1,1]. See [16] for more details. But, there is a special sub-case for the differential equation (2) whose solution is a basic class of symmetric orthogonal polynomials having an *explicit form* together with four free parameters. In this section, we are going to introduce this class and obtain all its standard properties. Of course, because of rather tedious calculations, we withdraw the extra computations and only express the original results.

## 2. A Basic Class of Symmetric Orthogonal Polynomials (BCSOP) using the extended Sturm-Liouville theorem

By referring to the generic differential equation (2) let us consider the following options

$$A(x) = x^2(px^2+q) \ ; \quad An\ even\ function,$$
$$B(x) = x(rx^2+s) \quad ; \quad An\ odd\ function,$$
$$C(x) = x^2 \qquad\qquad ; \quad An\ even\ function, \qquad (9)$$
$$D(x) = 0 \qquad\qquad ; \quad An\ even\ function,$$
$$E(x) = -s \qquad\qquad ; \quad An\ even\ function,$$

where $p, q, r$ and $s$ are real free parameters and $\lambda_n = -n(r + (n-1)p)$. Therefore, we deal with a second order differential equation in the form

$$x^2(px^2+q)\Phi_n''(x) + x(rx^2+s)\Phi_n'(x) - \left(n(r+(n-1)p)x^2 + (1-(-1)^n)s/2\right)\Phi_n(x) = 0. \quad (10)$$

To derive the polynomial solution of this equation, first suppose $n = 2m$, i.e.

$$x(px^2+q)\Phi_{2m}''(x) + (rx^2+s)\Phi_{2m}'(x) - 2m(r+(2m-1)p)\,x\,\Phi_{2m}(x) = 0. \quad (11)$$

Solving this equation eventually yields

$$S_{2m}\!\left(\begin{matrix} r & s \\ p & q \end{matrix}\bigg| x\right) = \sum_{k=0}^{m}\binom{m}{k}\left(\prod_{j=0}^{m-(k+1)}\frac{(2j-1+2m)p+r}{(2j+1)q+s}\right)x^{2m-2k}, \quad (12)$$

where $\prod_{r=0}^{r=-1} a_r = 1$. Similarly, if $n = 2m+1$ is taken in (10), the simplified equation

$$x^2(px^2+q)\Phi_{2m+1}''(x) + x(rx^2+s)\Phi_{2m+1}'(x) - ((2m+1)(r+2mp)x^2 + s)\Phi_{2m+1}(x) = 0, \quad (13)$$

has the polynomial solution

$$S_{2m+1}\begin{pmatrix} r & s \\ p & q \end{pmatrix} x = \sum_{k=0}^{m} \binom{m}{k} \left( \prod_{j=0}^{m-(k+1)} \frac{(2j+1+2m)p+r}{(2j+3)q+s} \right) x^{2m+1-2k} . \tag{14}$$

Hence, combining (12) with (14) finally gives

$$S_n\begin{pmatrix} r & s \\ p & q \end{pmatrix} x = \sum_{k=0}^{[n/2]} \binom{[n/2]}{k} \left( \prod_{i=0}^{[n/2]-(k+1)} \frac{(2i+(-1)^{n+1}+2[n/2])p+r}{(2i+(-1)^{n+1}+2)q+s} \right) x^{n-2k} \tag{15}$$

which is the most common source of classical *symmetric* orthogonal polynomials with four free parameters $p, q, r$ and $s$ where neither both values $q$ and $s$ nor both values $p$ and $r$ can vanish together. Here we mention that almost all known symmetric orthogonal polynomials, such as Legendre polynomials, first and second kind of Chebyshev polynomials, ultraspherical polynomials, generalized ultraspherical polynomials (GUP), Hermite polynomials and generalized Hermite polynomials (GHP) are special sub-cases of (15) and can be expressed in terms of it directly. Because of this matter, let us call these polynomials, "The second kind of classical orthogonal polynomials". Furthermore, there are two other symmetric sequences of *finite* orthogonal polynomials that are special sub-cases of general representation (15) and will be introduced in the sections 5 and 6.

The symbol $S_n\begin{pmatrix} r & s \\ p & q \end{pmatrix} x$ is used only because of Symmetry property, however, to abbreviate this notation in the text, we will show it as $S_n(p,q,r,s;x)$. A straightforward result from (15) is that

$$S_{2n+1}\begin{pmatrix} r & s \\ p & q \end{pmatrix} x = x S_{2n}\begin{pmatrix} r+2p & s+2q \\ p & q \end{pmatrix} x . \tag{16}$$

Moreover, since the monic type of orthogonal polynomials (i.e. with leading coefficient 1) is often required, we define

$$\bar{S}_n\begin{pmatrix} r & s \\ p & q \end{pmatrix} x = \prod_{i=0}^{[n/2]-1} \frac{(2i+(-1)^{n+1}+2)q+s}{(2i+(-1)^{n+1}+2[n/2])p+r} S_n\begin{pmatrix} r & s \\ p & q \end{pmatrix} x . \tag{17}$$

For instance, if $n = 0, ..., 5$ in (17) we have

$$\bar{S}_0\begin{pmatrix} r & s \\ p & q \end{pmatrix} x = 1, \quad \bar{S}_1\begin{pmatrix} r & s \\ p & q \end{pmatrix} x = x, \quad \bar{S}_2\begin{pmatrix} r & s \\ p & q \end{pmatrix} x = x^2 + \frac{q+s}{p+r},$$

$$\bar{S}_3\begin{pmatrix} r & s \\ p & q \end{pmatrix} x = x^3 + \frac{3q+s}{3p+r} x,$$

$$\bar{S}_4\begin{pmatrix} r & s \\ p & q \end{pmatrix} x = x^4 + 2\frac{3q+s}{5p+r} x^2 + \frac{(3q+s)(q+s)}{(5p+r)(3p+r)},$$

$$\bar{S}_5\begin{pmatrix} r & s \\ p & q \end{pmatrix} x = x^5 + 2\frac{5q+s}{7p+r} x^3 + \frac{(5q+s)(3q+s)}{(7p+r)(5p+r)} x.$$

(17.1)

The explicit form of polynomials (17) helps us obtain a three-term recurrence relation for the mentioned polynomials. In other words, if we assume that (17) satisfies the following relation

$$\bar{S}_{n+1}(x) = x\bar{S}_n(x) + C_n\begin{pmatrix} r & s \\ p & q \end{pmatrix}\bar{S}_{n-1}(x) \quad ; \quad \bar{S}_0(x) = 1, \quad \bar{S}_1(x) = x, \quad n \in N, \qquad (18)$$

then after doing a series of computations in hand, we obtain

$$C_n\begin{pmatrix} r & s \\ p & q \end{pmatrix} = \frac{pqn^2 + \big((r-2p)q - (-1)^n ps\big)n + (r-2p)s(1-(-1)^n)/2}{(2pn+r-p)(2pn+r-3p)}, \qquad (18.1)$$

which in fact reveals the explicit form of (18).
Now, since the recurrence relation (18) is explicitly known, to determine the norm square value of the polynomials we can apply Favard's theorem [6] by noting that there is orthogonality with respect to a weight function. According to this theorem, if $\{P_n(x)\}_{n=0}^{\infty}$ is defined by

$$xP_n(x) = A_n P_{n+1}(x) + B_n P_n(x) + C_n P_{n-1}(x) \quad ; \quad n = 0,1,2,\ldots, \qquad (19)$$

where $P_{-1}(x) = 0$, $P_0(x) = 1$, $A_n, B_n, C_n$ real and $A_n C_{n+1} > 0$ for $n = 0,1,\ldots$, then there exists a weight function $W(x)$ so that

$$\int_{-\infty}^{\infty} W(x) P_n(x) P_m(x) \, dx = \left(\prod_{i=0}^{n-1} \frac{C_{i+1}}{A_i} \int_{-\infty}^{\infty} W(x) \, dx\right) \delta_{n,m}. \qquad (20)$$

Moreover, if the positive condition $A_n C_{n+1} > 0$ only holds for $n = 0,1,\ldots,N$ then the orthogonality relation (20) only holds for a finite number of $m, n$. This latter note will help us in the next sections to obtain two new sub-classes of $\bar{S}_n(p,q,r,s;x)$ that are *finitely* orthogonal on the interval $(-\infty, \infty)$. It is clear that the Favard theorem is also valid for the recurrence equation (18) in which $A_n = 1$, $B_n = 0$ and $C_n = -C_n(p,q,r,s)$. Furthermore, the condition $-C_{n+1}(p,q,r,s) > 0$ must always be satisfied if one demands to apply the Favard theorem for (18). By noting these comments and (20), the generic form of the orthogonality relation of BCSOP can be designed as

(21)
$$\int_{-\alpha}^{\alpha} W\begin{pmatrix} r & s \\ p & q \end{pmatrix} x\Big) \bar{S}_n\begin{pmatrix} r & s \\ p & q \end{pmatrix} x\Big) \bar{S}_m\begin{pmatrix} r & s \\ p & q \end{pmatrix} x\Big) dx = \left((-1)^n \prod_{i=1}^{i=n} C_i\begin{pmatrix} r & s \\ p & q \end{pmatrix} \int_{-\alpha}^{\alpha} W\begin{pmatrix} r & s \\ p & q \end{pmatrix} x\Big) dx\right) \delta_{n,m},$$

where the weight function, by referring to relations (4) and (9), is defined as

$$W\begin{pmatrix} r & s \\ p & q \end{pmatrix} x\Big) = x^2 \exp\left(\int \frac{(r-4p)x^2 + (s-2q)}{x(px^2+q)} dx\right) = \exp\left(\int \frac{(r-2p)x^2 + s}{x(px^2+q)} dx\right), \qquad (21.1)$$

and $\alpha$ takes the standard values $1, \infty$. Note that the function $(px^2 + q)W(p,q,r,s;x)$ must vanish at $x = \alpha$ in order to establish the main orthogonality relation (21).

## 2.1. An analogue of Pearson distributions family

The positive function (21.1) can be investigated from statistical point of view too. In fact, this function is an analogue of Pearson distributions family having the general form

$$\rho\begin{pmatrix} d & e \\ a & b & c \end{pmatrix} x = \exp(\int \frac{(d-2a)x + (e-b)}{ax^2 + bx + c} dx), \tag{22}$$

and satisfying the first order differential equation

$$\frac{d}{dx}((ax^2 + bx + c)\rho(x)) = (dx + e)\rho(x). \tag{22.1}$$

Hence, we would like to point out that, similar to equation (22.1), the weight function (21.1) satisfies a first order differential equation as

$$x \frac{d}{dx}((px^2 + q)W(x)) = (rx^2 + s)W(x), \tag{23}$$

which is equivalent to

$$\frac{d}{dx}(x^2(px^2 + q)W(x)) = x(r^*x^2 + s^*)W(x) \quad \text{s.t.} \quad \begin{cases} r^* = r + 2p \\ s^* = s + 2q \end{cases}. \tag{23.1}$$

From relations (23) or (23.1) one can deduce that $W(p,q,r,s;x)$ is an analytical integrable function and since it is also a positive function, its probability density function (pdf) must be available.

In general, there are four main sub-classes of distributions family (21.1) ( and consequently sub-solutions of equation (23)) whose explicit pdfs are respectively as follows

$$K_1 W\begin{pmatrix} -2a-2b-2, & 2a \\ -1, & 1 \end{pmatrix} x = \frac{\Gamma(a+b+3/2)}{\Gamma(a+1/2)\Gamma(b+1)} x^{2a}(1-x^2)^b \ ; \ -1 \le x \le 1; \tag{24}$$

$$a + 1/2 > 0 \ ; \ b + 1 > 0.$$

$$K_2 W\begin{pmatrix} -2, & 2a \\ 0, & 1 \end{pmatrix} x = \frac{1}{\Gamma(a+1/2)} x^{2a} e^{-x^2} \ ; \ -\infty < x < \infty \ ; \ a + 1/2 > 0. \tag{25}$$

$$K_3 W\begin{pmatrix} -2a-2b+2, & -2a \\ 1, & 1 \end{pmatrix} x = \frac{\Gamma(b)}{\Gamma(b+a-1/2)\Gamma(-a+1/2)} \frac{x^{-2a}}{(1+x^2)^b} \ ; \ -\infty < x < \infty \ ; \tag{26}$$

$$b + a > 1/2 \ ; \ a < 1/2 \ ; \ b > 0.$$

$$K_4 W\begin{pmatrix} -2a+2, & 2 \\ 1, & 0 \end{pmatrix} x = \frac{1}{\Gamma(a-1/2)} x^{-2a} e^{-\frac{1}{x^2}} \ ; \ -\infty < x < \infty \ ; \ a > 1/2. \tag{27}$$

The values $K_i$ ; $i = 1,2,3,4$ play the normalizing constant role in above distributions. Moreover, as it is observed, the value of distribution vanishes at $x = 0$ in each four cases, i.e. $W(p,q,r,s;0) = 0$ for $s \neq 0$. Because of this property, we call (21.1) "The dual symmetric distributions family".

**2.2. A generalization of dual symmetric distributions family**
First it is important to note that if $s = 0$ in (23), the foresaid equation will be reduced to a special sub-case of Pearson differential equation (22.1). Hence, we hereafter suppose that $s \neq 0$. Since the explicit forms of $S_n(p,q,r,s;x)$ in (15), $C_n(p,q,r,s)$ in (18.1) and $W(p,q,r,s;x)$ in (21.1) are all known, a further basic pdf can be defined by referring to the orthogonality relation (21) directly, so that we have

$$D_m\!\left(\begin{matrix} r & s \\ p & q \end{matrix}\bigg|x\right) = \frac{(-1)^m}{\prod_{i=1}^{i=m} C_i\!\left(\begin{matrix} r & s \\ p & q \end{matrix}\right) \int_{-\alpha}^{\alpha} W\!\left(\begin{matrix} r & s \\ p & q \end{matrix}\bigg|x\right) dx} \; W\!\left(\begin{matrix} r & s \\ p & q \end{matrix}\bigg|x\right) \times \left(\bar{S}_m\!\left(\begin{matrix} r & s \\ p & q \end{matrix}\bigg|x\right)\right)^2. \tag{28}$$

Clearly choosing $m = 0$ in this definition gives the same as dual symmetric distributions family. Here let us add that the Fisher information [9] and Boltzmann-Shannon information entropy [7], as two important factors in statistical estimation theory, can be studied for the generalized distribution (28).

**2.3. A direct relationship between first and second kind of classical orthogonal polynomials**
One can verify that there is a direct relation between first kind of classical orthogonal polynomials (including Jacobi, Laguerre and Hermite polynomials as well as three finite classes of orthogonal polynomials, see e.g. [14,15] for more details) and the explicit polynomials $S_n(p,q,r,s;x)$ indicated in (15). To study the first kind of classical orthogonal polynomials we refer the readers to a valuable book written by Nikiforov and Uvarov [19].
To find this relationship, we start with the following differential equation

$$(ax^2 + bx + c)y_n''(x) + (dx + e)y_n'(x) - n((n-1)a + d)y_n(x) = 0 \;. \tag{29}$$

According to [19], the monic polynomial solution of the equation (29) can be shown by a Rodrigues-type formula as

$$\bar{P}_n\!\left(\begin{matrix} d & e \\ a & b & c \end{matrix}\bigg|x\right) = \frac{1}{(\prod_{k=1}^{n} d + (n+k-2)a)\rho\!\left(\begin{matrix} d & e \\ a & b & c \end{matrix}\bigg|x\right)} \times \frac{d^n}{dx^n}((ax^2+bx+c)^n \rho\!\left(\begin{matrix} d & e \\ a & b & c \end{matrix}\bigg|x\right)) \;, \tag{30}$$

where $\rho(a,b,c,d,e;x)$ is the same definition as (22).

On the other hand, since $\bar{S}_{2n}(p,q,r,s;x)$ is generally an even function, taking $x = wt^2 + v$ in (29) gives (31)

$$t^2(aw^2t^4 + w(2av+b)t^2 + av^2 + bv + c)y''_n(t) + t((2d-a)w^2t^4 + (2wv(d-a) + w(2e-b))t^2$$
$$- (av^2 + bv + c))y'_n(t) - 4w^2 n(d + (n-1)a)t^4 y_n(t) = 0.$$

If (31) is equated with (11), we should have

$$av^2 + bv + c = 0 \quad or \quad v = \frac{-b \pm \sqrt{b^2 - 4ac}}{2a}. \tag{31.1}$$

The condition (31.1) simplifies the equation (31) to

$$t(awt^2 \pm \sqrt{b^2 - 4ac})y''_n(t) + ((2d-a)wt^2 + (\frac{d}{a}-1)(-b \pm \sqrt{b^2 - 4ac}) + 2e - b)y'_n(t)$$
$$- 4wn(d + (n-1)a)t\, y_n(t) = 0 \Leftrightarrow y_n(t) = \bar{P}_n\left(\begin{array}{cc} d & e \\ a & b & c \end{array}\middle| wt^2 + \frac{-b \pm \sqrt{b^2 - 4ac}}{2a}\right). \tag{31.2}$$

The equation (31.2) is clearly a special case of (11). This means that

(32)
$$\bar{P}_n\left(\begin{array}{cc} d & e \\ a & b & c \end{array}\middle| wt^2 + \frac{-b \pm \sqrt{b^2 - 4ac}}{2a}\right) = K\, \bar{S}_{2n}\left(\begin{array}{cc} (2d-a)w, & (\frac{d}{a}-1)(-b \pm \sqrt{b^2-4ac}) + 2e - b \\ aw, & \pm\sqrt{b^2 - 4ac} \end{array}\middle| t\right),$$

where $K$ is the leading coefficient of the left-hand side polynomial of relation (32) divided to leading coefficient of the right-hand side polynomial.
But as we observe, there exist 5 free parameters $a, b, c, d, e$ in the left-hand side of (32). So, one of them must be pre-assigned in order that one can get the explicit form of $\bar{S}_{2n}(p,q,r,s;x)$ in terms of $\bar{P}_n(a,b,c,d,e;wt^2+v)$ similarly. For this purpose, if for instance $c = 0$ is considered in (32), the two following cases appear

$$\bar{S}_{2n}\left(\begin{array}{cc} r & s \\ p & q \end{array}\middle| t\right) = w^{-n} \bar{P}_n\left(\begin{array}{cc} \frac{r+p}{2w}, & \frac{s+q}{2} \\ p/w, & q, & 0 \end{array}\middle| wt^2\right) \quad ; \quad w \neq 0, \tag{33}$$

$$\bar{S}_{2n}\left(\begin{array}{cc} r & s \\ p & q \end{array}\middle| t\right) = w^{-n} \bar{P}_n\left(\begin{array}{cc} \frac{r+p}{2w}, & \frac{sp-rq}{2p} \\ p/w, & -q, & 0 \end{array}\middle| wt^2 + w\frac{q}{p}\right) \quad ; \quad p, w \neq 0. \tag{34}$$

Furthermore, if (16) is applied for two latter relations we respectively get

$$\bar{S}_{2n+1}\begin{pmatrix} r & s \\ p & q \end{pmatrix} t = w^{-n} t \, \bar{P}_n \begin{pmatrix} \dfrac{r+3p}{2w}, & \dfrac{s+3q}{2} \\ p/w, & q, & 0 \end{pmatrix} wt^2 \quad ; \quad w \neq 0 , \qquad (35)$$

$$\bar{S}_{2n+1}\begin{pmatrix} r & s \\ p & q \end{pmatrix} t = w^{-n} t \, \bar{P}_n \begin{pmatrix} \dfrac{r+3p}{2w}, & \dfrac{sp-rq}{2p} \\ p/w, & -q, & 0 \end{pmatrix} wt^2 + w\dfrac{q}{p} \quad ; \quad p, w \neq 0 . \qquad (36)$$

### 2.4. Some further standard properties of BCSOP

The relations (33)-(36) are useful tool to get a generating function for BCSOP. Usually, a generating function for a system of polynomials $P_n(z)$ is defined by a function like $G(z,t)$ whose expansion in powers of $t$ has, for sufficiently small $|t|$, the form

$$G(z,t) = \sum_{n=0}^{\infty} P_n(z) \frac{t^n}{n!} . \qquad (37)$$

If $P_n(z)$ has the Rodrigues-type formula [1,19], Cauchy's integral formula

$$f^{(n)}(z) = \frac{n!}{2\pi i} \int_{\Omega} \frac{f(u)}{(u-z)^{n+1}} du , \qquad (38)$$

, where $\Omega$ is a closed contour surrounding the point $u = z$, is employed to obtain $G(z,t)$ explicitly, see e.g. [19, p. 27]. This means that by considering the Rodrigues-type representation (30) for $c = 0$ and applying Cauchy's integral theorem on it we have

$$\sum_{n=0}^{\infty} (\prod_{k=1}^{n} d + (n+k-2)a) \bar{P}_n \begin{pmatrix} d & e \\ a & b & 0 \end{pmatrix} z \frac{t^n}{n!} = \frac{\rho\begin{pmatrix} d & e \\ a, & b, & 0 \end{pmatrix} u}{\rho\begin{pmatrix} d & e \\ a, & b, & 0 \end{pmatrix} z} \times \frac{1}{\sqrt{(1-bt)^2 - 4atz}} . \qquad (39)$$

$$where: \quad u = \frac{1 - bt + \sqrt{(1-bt)^2 - 4atz}}{2at} .$$

If $z \to z^2$, $t \to t^2$, $a = p$, $b = q$, $d = \dfrac{r+p}{2}$ and $e = \dfrac{s+q}{2}$ are substituted into (39), then by noting the relations (33)-(36) we have $\qquad (40)$

$$\sum_{n=0}^{\infty} (\prod_{k=1}^{n} \frac{r+p}{2} + (n+k-2)p) \bar{S}_{2n} \begin{pmatrix} r & s \\ p & q \end{pmatrix} z \frac{t^{2n}}{n!} = \frac{\rho\begin{pmatrix} (r+p)/2 & (s+q)/2 \\ p, & q, & 0 \end{pmatrix} u^*}{\rho\begin{pmatrix} (r+p)/2 & (s+q)/2 \\ p, & q, & 0 \end{pmatrix} z^2} \frac{1}{\sqrt{(1-qt^2)^2 - 4pt^2 z^2}}$$

$$where: \quad u^* = \frac{1 - qt^2 + \sqrt{(1-qt^2)^2 - 4pt^2 z^2}}{2pt^2} .$$

Similarly by multiplying both sides of (39) by $tz$ and using (35) we get

$$\sum_{n=0}^{\infty}(\prod_{k=1}^{n}\frac{r+3p}{2}+(n+k-2)p)\bar{S}_{2n+1}\left(\begin{array}{cc} r & s \\ p & q \end{array}\bigg| z\right)\frac{t^{2n+1}}{n!} = \frac{\rho\left(\begin{array}{cc} (r+3p)/2 & (s+3q)/2 \\ p, \ q, \ 0 \end{array}\bigg| u^*\right)}{\rho\left(\begin{array}{cc} (r+3p)/2 & (s+3q)/2 \\ p, \ q, \ 0 \end{array}\bigg| z^2\right)}\cdot\frac{tz}{\sqrt{(1-qt^2)^2-4pt^2z^2}}. \tag{41}$$

Now if we define

$$S_n^*\left(\begin{array}{cc} r & s \\ p & q \end{array}\bigg| x\right) = \frac{\prod_{k=1}^{[n/2]} r/2 + ((2n-5+(-1)^{n+1})/4+k)p}{[n/2]!}\bar{S}_n\left(\begin{array}{cc} r & s \\ p & q \end{array}\bigg| x\right), \tag{42}$$

then by noting (40) and (41) a kind of generating function for BCSOP can finally be derived as follows

$$\sum_{n=0}^{\infty}S_n^*\left(\begin{array}{cc} r & s \\ p & q \end{array}\bigg| z\right)t^n = \frac{1}{\sqrt{(1-qt^2)^2-4pt^2z^2}}\left(\frac{\rho\left(\begin{array}{cc} (r+p)/2, & (s+q)/2 \\ p, \ q, \ 0 \end{array}\bigg| u^*\right)}{\rho\left(\begin{array}{cc} (r+p)/2, & (s+q)/2 \\ p, \ q, \ 0 \end{array}\bigg| z^2\right)} + \frac{(tz)\rho\left(\begin{array}{cc} (r+3p)/2, & (s+3q)/2 \\ p, \ q, \ 0 \end{array}\bigg| u^*\right)}{\rho\left(\begin{array}{cc} (r+3p)/2, & (s+3q)/2 \\ p, \ q, \ 0 \end{array}\bigg| z^2\right)}\right) \tag{43}$$

where $u^*$ is the same form as (40).

The explicit form of polynomials (15) can also be applied to show a generic hypergeometric representation, because by using it one can easily indicate the coefficients of polynomials (15) in terms of the Pochhammer symbol: $(a)_m = a(a+1)...(a+m-1)$. For instance, we have

$$\prod_{i=0}^{[n/2]-(k+1)}(2i+(-1)^{n+1}+2[n/2])p+r = (p/2)^{[n/2]-k}\left(r/2p+[n/2]+(-1)^{n+1}/2\right)_{[\frac{n}{2}]-k}. \tag{44}$$

Hence, after doing some calculations, we eventually find that

$$\bar{S}_n\left(\begin{array}{cc} r & s \\ p & q \end{array}\bigg| x\right) = x^n {}_2F_1\left(\begin{array}{c} -[n/2], \ (q-s)/2q-[(n+1)/2] \\ -(r+(2n-3)p)/2p \end{array}\bigg| -\frac{q}{px^2}\right), \tag{45}$$

Where ${}_2F_1\left(\begin{array}{cc} \alpha & \beta \\ & \gamma \end{array}\bigg| x\right) = \sum_{k=0}^{\infty}\frac{(\alpha)_k(\beta)_k}{(\gamma)_k}\frac{x^k}{k!}$ denotes the hypergeometric function of order (2,1) [12].

On the other hand, since

$${}_2F_1\left(\begin{array}{cc} a & b \\ & c \end{array}\bigg| z\right) = \frac{\Gamma(c)}{\Gamma(b)\Gamma(c-b)}\int_0^1 t^{b-1}(1-t)^{c-b-1}(1-tz)^{-a}dt \Leftrightarrow \text{Re}\,c > \text{Re}\,b > 0\ ;\ |z|<1, \tag{46}$$

the integral representation of BCSOP will easily be derived by referring to (45). In this way, by applying the Gauss identity [12]

$${}_2F_1\left(\begin{array}{cc} a & b \\ & c \end{array}\bigg| 1\right) = \frac{\Gamma(c)\Gamma(c-a-b)}{\Gamma(c-a)\Gamma(c-b)}, \tag{47}$$

one can also determine the value of polynomials $\bar{S}_n(p,q,r,s;x)$ at a specific boundary point, i.e. $x = \sqrt{-q/p}$. To reach this goal, it is sufficient to put $x = \sqrt{-q/p}$ in (45) and use (47) to get

$$\bar{S}_n\left(\begin{array}{cc} r & s \\ p & q \end{array}\middle|\sqrt{-\frac{q}{p}}\right) = \frac{(-\frac{q}{p})^{\frac{n}{2}} \Gamma(\frac{3}{2} - \frac{r}{2p})\Gamma(1 + \frac{s}{2q} - \frac{r}{2p})}{(\frac{3}{2} - \frac{r}{2p} - n)_n \Gamma(\frac{3}{2} - \frac{r}{2p} + [\frac{n}{2}] - n) \Gamma(1 + \frac{s}{2q} - \frac{r}{2p} - [\frac{n}{2}])} . \tag{48}$$

Note that to derive the above identity, the following equalities have been used

$$[\frac{n+1}{2}] - [\frac{n}{2}] = \frac{1-(-1)^n}{2}, \qquad [\frac{n+1}{2}] + [\frac{n}{2}] = n,$$

$$[\frac{n}{2}] = \frac{n}{2} - \frac{1-(-1)^n}{4}, \qquad (a)_n = (-1)^n (1-a-n)_n. \tag{49}$$

Now the standard properties of BCSOP have been found and it is a good position to introduce four special sub-cases of basic polynomials (15) in detail.

## 3. First subclass, Generalized ultraspherical polynomials (GUP)

These polynomials were first investigated by Chihara in detail [see 6]. He obtained the main properties of GUP via a direct relation between them and Jacobi orthogonal polynomials. The asymptotic behaviors of foresaid polynomials were also studied by Konoplev [13]. On the other hand, Charris and Ismail in [5] (see also [11]) introduced the Sieved random walk polynomials to show that the generalized Ultraspherical polynomials are a special case of them. Of course, there are some other generalizations of Ultraspherical polynomials. For instance, Askey in [4] introduced two classes of orthogonal polynomials as a limiting case of the q-Wilson polynomials on $[-1,1]$ with the weight functions

$$W_1(x,\alpha,\lambda) = (1-x^2)^{\alpha-\frac{1}{2}} |U_{k-1}(x)|^{2\alpha} |T_k(x)|^\lambda$$

and

$$W_2(x,\alpha,\lambda) = (1-x^2)^{\alpha+\frac{1}{2}} |U_{k-1}(x)|^{2\alpha} |T_k(x)|^\lambda$$

in which $k \in N$ is a fixed integer and $T_k(x)$ and $U_{k-1}(x)$ are respectively the Chebyshev polynomials of the first and second kind, to generalize GUP for $k=1$ and $\lambda=0$. For the case $\lambda=0$ the polynomials were introduced by Al-Salam, Allaway and Askey [2] as a limiting case of the q-Ultraspherical polynomials of Rogers [21].

Now, in this section, we intend to show that GUP can directly be represented in terms of $\bar{S}_n(p,q,r,s;x)$ and consequently all its standard properties can be obtained. For this purpose, it is only enough to have the initial vector corresponding to these polynomials and replace it into the standard properties of BCSOP.

**Definition**

Choose the initial vector $(p,q,r,s) = (-1,1,-2a-2b-2,2a)$ and substitute it into (15) to get

$$\bar{S}_n\left(\begin{matrix}-2a-2b-2, & 2a \\ -1, & 1\end{matrix}\bigg| x\right) = \prod_{i=0}^{[n/2]-1} \frac{2i+2a+2-(-1)^n}{-2i-(2b+2a+2-(-1)^n+2[n/2])}$$
$$\times \sum_{k=0}^{[n/2]} \binom{[n/2]}{k} \left(\prod_{i=0}^{[\frac{n}{2}]-(k+1)} \frac{-2i-(2b+2a+2-(-1)^n+2[n/2])}{2i+2a+2-(-1)^n}\right) x^{n-2k}, \quad (50)$$

as the explicit form of monic GUP. Furthermore, since the Ultraspherical (Gegenbauer), Legendre, and Chebyshev polynomials of the first and second kind are all special cases of GUP, they can be denoted by $\bar{S}_n(p,q,r,s;x)$ directly and we have

Ultraspherical polynomials:

$$C_n^a(x) = \frac{2^n (a)_n}{n!} \bar{S}_n\left(\begin{matrix}-2a-1 & 0 \\ -1 & 1\end{matrix}\bigg| x\right). \quad (50.1)$$

Legendre polynomials:

$$P_n(x) = \frac{(2n)!}{(n!)^2 2^n} \bar{S}_n\left(\begin{matrix}-2 & 0 \\ -1 & 1\end{matrix}\bigg| x\right). \quad (50.2)$$

Chebyshev polynomials of the first kind:

$$T_n(x) = 2^{n-1} \bar{S}_n\left(\begin{matrix}-1 & 0 \\ -1 & 1\end{matrix}\bigg| x\right). \quad (50.3)$$

Chebyshev polynomials of the first kind:

$$U_n(x) = 2^n \bar{S}_n\left(\begin{matrix}-3 & 0 \\ -1 & 1\end{matrix}\bigg| x\right). \quad (50.4)$$

**Recurrence relation of monic polynomials**

By replacing the initial vector (50) into the explicit expression $C_n(p,q,r,s)$, given in (18.1), the recurrence relation of monic GUP takes the form (51)

$$\bar{S}_{n+1}(x) = x\bar{S}_n(x) + C_n\left(\begin{matrix}-2a-2b-2, & 2a \\ -1, & 1\end{matrix}\right) \bar{S}_{n-1}(x) \quad ; \quad \bar{S}_0(x) = 1, \; \bar{S}_1(x) = x, \; n \in N,$$

where according to (18.1)

$$C_n\left(\begin{matrix}-2a-2b-2, & 2a \\ -1, & 1\end{matrix}\right) = \frac{-n^2 - (2b+2(1-(-1)^n)a)n - 2a(a+b)(1-(-1)^n)}{(2n+2a+2b-1)(2n+2a+2b+1)}$$
$$= \frac{-(n+(1-(-1)^n)a)(n+(1-(-1)^n)a+2b)}{(2n+2a+2b-1)(2n+2a+2b+1)}. \quad (51.1)$$

**Orthogonality relation**

Clearly the weight function of GUP is the same distribution as (24) without considering its normalizing constant, i.e. $x^{2a}(1-x^2)^b$. Also, since this function must be even and positive, the condition $(-1)^{2a}=1$ is essential. Hence, the mentioned weight function can also be considered as $|x|^{2a}(1-x^2)^b$; $x \in [-1,1]$. By noting this comment and generic relation (21) for $\alpha = 1$, the orthogonality relation of first sub-class takes the form

$$\int_{-1}^{1} x^{2a}(1-x^2)^b \bar{S}_n\left(\begin{array}{cc} -2a-2b-2, & 2a \\ -1, & 1 \end{array}\bigg| x\right) \bar{S}_m\left(\begin{array}{cc} -2a-2b-2, & 2a \\ -1, & 1 \end{array}\bigg| x\right) dx$$
$$= \left((-1)^n \prod_{i=1}^{i=n} C_i\left(\begin{array}{cc} -2a-2b-2, & 2a \\ -1, & 1 \end{array}\right) \int_{-1}^{1} x^{2a}(1-x^2)^b dx\right) \delta_{n,m}, \quad (52)$$

where

$$\int_{-1}^{1} x^{2a}(1-x^2)^b dx = B(a+\frac{1}{2}, b+1) = \frac{\Gamma(a+1/2)\Gamma(b+1)}{\Gamma(a+b+3/2)}. \quad (52.1)$$

From (52.1) one can conclude that the constraints on the parameters $a$ and $b$ should be $a+1/2 > 0$, $(-1)^{2a}=1$ and $b+1 > 0$. Note that $B(\lambda_1, \lambda_2)$ in (52.1) denotes the Beta integral [3] having various definitions as

$$B(\lambda_1, \lambda_2) = \int_0^1 x^{\lambda_1-1}(1-x)^{\lambda_2-1} dx = \int_{-1}^1 x^{2\lambda_1-1}(1-x^2)^{\lambda_2-1} dx = \int_0^\infty \frac{x^{\lambda_1-1}}{(1+x)^{\lambda_1+\lambda_2}} dx$$
$$= 2\int_0^{\pi/2} Sin^{(2\lambda_1-1)}x Cos^{(2\lambda_2-1)}x\, dx = \frac{\Gamma(\lambda_1)\Gamma(\lambda_2)}{\Gamma(\lambda_1+\lambda_2)} = B(\lambda_2, \lambda_1). \quad (53)$$

**Differential equation:** $\Phi_n(x) = S_n(-1,1,-2a-2b-2,2a;x)$.

To derive the differential equation of GUP it is enough to substitute the initial vector (50) into the main differential equation (10) to get to (54)

$$x^2(-x^2+1)\Phi_n''(x) - 2x((a+b+1)x^2 - a)\Phi_n'(x) + \left(n(2a+2b+n+1)x^2 + ((-1)^n-1)a\right)\Phi_n(x) = 0.$$

### 3.1. Fifth and Sixth kind of Chebyshev polynomials

As we know, four kinds of trigonometric orthogonal polynomials, known as first, second, third and fourth kind of Chebyshev polynomials, have been investigated in the literature up to now, see e.g. [20,10,6,22]. The explicit definitions of these polynomials are respectively

$$T_n(x) = 2^{n-1}\prod_{k=1}^n (x - Cos\frac{(2k-1)\pi}{2n}) = Cos(n\theta) \quad ; \quad x = Cos\theta, \quad (55)$$

$$U_n(x) = 2^n \prod_{k=1}^n (x - Cos\frac{k\pi}{n+1}) = \frac{Sin((n+1)\theta)}{Sin\theta} \quad ; \quad x = Cos\theta, \quad (56)$$

$$V_n(x) = 2^n \prod_{k=1}^{n}(x - Cos\frac{(2k-1)\pi}{2n+1}) = \frac{Cos((2n+1)\theta)}{Cos(\theta)} \quad ; \quad x = Cos(2\theta) , \tag{57}$$

$$W_n(x) = 2^n \prod_{k=1}^{n}(x - Cos\frac{2k\pi}{2n+1}) = \frac{Sin((2n+1)\theta)}{Sin(\theta)} \quad ; \quad x = Cos(2\theta) . \tag{58}$$

But, in this section we would like to add that there exist two further kinds of Half-trigonometric orthogonal polynomials, which are particular sub-cases of $\overline{S}_n(p,q,r,s;x)$. Since they are generated by using the first and second kind of Chebyshev polynomials and have the half-trigonometric forms, let us call them the fifth and sixth kind of Chebyshev polynomials.
To generate these two sequences, we should refer to the important relation (16). According to (50.3), the initial vector of first kind Chebyshev polynomials is: $(p,q,r,s) = (-1,1,-1,0)$. So, if this vector is replaced into (16) then

$$\overline{S}_{2n+1}\left(\begin{matrix} -1 & 0 \\ -1 & 1 \end{matrix} \middle| x\right) = \overline{T}_{2n+1}(x) = x\,\overline{S}_{2n}\left(\begin{matrix} -3 & 2 \\ -1 & 1 \end{matrix} \middle| x\right). \tag{59}$$

By means of (59), the secondary vector $(p,q,r,s) = (-1,1,-3,2)$, as a special case of the set of vectors $(p,q,r,s)$, appears. Using this new vector first we define the half-trigonometric polynomials

$$\overline{X}_n(x) = \overline{S}_n\left(\begin{matrix} -3 & 2 \\ -1 & 1 \end{matrix} \middle| x\right) = \begin{cases} \dfrac{(-1)^{n/2}}{n+1}\dfrac{Cos((n+1)\theta)}{Cos\theta} & \text{if } n = 2m, \\ \overline{S}_n(-1,1,-3,2;x) & \text{if } n = 2m+1, \end{cases} \tag{60}$$

where $Cos\theta = x$. According to (18) and (18.1), $\overline{X}_n(x)$ satisfies the recurrence relation

$$\overline{X}_{n+1}(x) = x\,\overline{X}_n(x) + C_n\left(\begin{matrix} -3 & 2 \\ -1 & 1 \end{matrix}\right)\overline{X}_{n-1}(x) \quad ; \quad \overline{X}_0(x) = 1 , \overline{X}_1(x) = x , n \in N, \tag{61}$$

in which

$$C_n\left(\begin{matrix} -3 & 2 \\ -1 & 1 \end{matrix}\right) = \frac{-(n-(-1)^n)(n+1-(-1)^n)}{4n(n+1)} . \tag{61.1}$$

Consequently, substituting the secondary vector (-1,1,-3,2) into the generic relation (21) gives the orthogonality relation of the first kind of half-trigonometric orthogonal polynomials as

$$\int_{-1}^{1} W\left(\begin{matrix} -3 & 2 \\ -1 & 1 \end{matrix} \middle| x\right)\overline{X}_n(x)\overline{X}_m(x)dx = \left((-1)^n \prod_{i=1}^{n} C_i\left(\begin{matrix} -3 & 2 \\ -1 & 1 \end{matrix}\right)\int_{-1}^{1} W\left(\begin{matrix} -3 & 2 \\ -1 & 1 \end{matrix} \middle| x\right)dx\right)\delta_{n,m} . \tag{62}$$

On the other hand since

$$\int_{-1}^{1} W\begin{pmatrix} -3 & 2 \\ -1 & 1 \end{pmatrix} x \, dx = \int_{-1}^{1} \frac{x^2}{\sqrt{1-x^2}} \, dx = B(\frac{3}{2}, \frac{1}{2}) = \frac{\pi}{2}, \quad (63)$$

(62) is simplified to

$$\int_{-1}^{1} \frac{x^2}{\sqrt{1-x^2}} \overline{X}_n(x) \overline{X}_m(x) \, dx = \left( (-1)^n \prod_{i=1}^{n} C_i \begin{pmatrix} -3 & 2 \\ -1 & 1 \end{pmatrix} \right) \frac{\pi}{2} \delta_{n,m}. \quad (64)$$

Clearly a half of polynomials $\overline{X}_n(x)$ is decomposable, i.e. we have

$$\overline{X}_{2n}(x) = \prod_{k=1}^{2n} (x - Cos\frac{(2k-1)\pi}{2(2n+1)}). \quad (65)$$

Therefore, if we can find the roots of $\overline{X}_{2n+1}(x)$ too, these polynomials will be useful in many numerical analysis applications such as Gaussian quadrature rules [6,22].
Similarly, the subject holds for the initial vector of the second kind Chebyshev polynomials, i.e. $(p,q,r,s) = (-1,1,-3,0)$. Again, if this vector is substituted into (16) then

$$\overline{S}_{2n+1}\begin{pmatrix} -3 & 0 \\ -1 & 1 \end{pmatrix} x = \overline{U}_{2n+1}(x) = x \overline{S}_{2n}\begin{pmatrix} -5 & 2 \\ -1 & 1 \end{pmatrix} x, \quad (66)$$

and subsequently the secondary vector is obtained in the form $(p,q,r,s) = (-1,1,-5,2)$. Now, by assuming that $Cos\theta = x$ let us define the polynomials

$$Y_n(x) = S_n\begin{pmatrix} -5 & 2 \\ -1 & 1 \end{pmatrix} x = \begin{cases} \frac{(-1)^{n/2}}{n+2} \frac{Sin((n+2)\theta)}{Cos\theta \, Sin\theta} & \text{if } n = 2m, \\ S_n(-1,1,-5,2;x) & \text{if } n = 2m+1, \end{cases} \quad (67)$$

satisfying the recurrence relation

$$\overline{Y}_{n+1}(x) = x \overline{Y}_n(x) + C_n\begin{pmatrix} -5 & 2 \\ -1 & 1 \end{pmatrix} \overline{Y}_{n-1}(x) \quad ; \quad \overline{Y}_0(x) = 1, \overline{Y}_1(x) = x, n \in N, \quad (68)$$

where

$$C_n\begin{pmatrix} -5 & 2 \\ -1 & 1 \end{pmatrix} = \frac{-(n+1-(-1)^n)(n+2-(-1)^n)}{4(n+1)(n+2)}, \quad (68.1)$$

and having the orthogonality relation

$$\int_{-1}^{1} W\begin{pmatrix} -5 & 2 \\ -1 & 1 \end{pmatrix} x \overline{Y}_n(x) \overline{Y}_m(x) \, dx = \left( (-1)^n \prod_{i=1}^{n} C_i \begin{pmatrix} -5 & 2 \\ -1 & 1 \end{pmatrix} \right) \int_{-1}^{1} W\begin{pmatrix} -5 & 2 \\ -1 & 1 \end{pmatrix} x \, dx \, \delta_{n,m}, \quad (69)$$

where

$$\int_{-1}^{1} W\begin{pmatrix} -5 & 2 \\ -1 & 1 \end{pmatrix} x \, dx = \int_{-1}^{1} x^2 \sqrt{1-x^2} \, dx = B(\frac{3}{2}, \frac{3}{2}) = \frac{\pi}{8} \, . \tag{70}$$

The relation (70) simplifies (69) as

$$\int_{-1}^{1} x^2 \sqrt{1-x^2} \, \overline{Y}_n(x) \overline{Y}_m(x) dx = \left( (-1)^n \prod_{i=1}^{n} C_i \begin{pmatrix} -5 & 2 \\ -1 & 1 \end{pmatrix} \right) \frac{\pi}{8} \delta_{n,m} \, . \tag{71}$$

Similar to previous case, $\overline{Y}_{2n}(x)$ is decomposable as

$$\overline{Y}_{2n}(x) = \prod_{k=1}^{2n} (x - Cos \frac{k\pi}{2n+2}) \, . \tag{72}$$

Here we should add that there are two other sequences of half-trigonometric polynomials that are *not orthogonal*, but can be shown in terms of BCSOP. These sequences are defined as

$$\overline{S}_n \begin{pmatrix} 1 & -2 \\ -1 & 1 \end{pmatrix} x = \begin{cases} \overline{S}_n(-1,1,1,-2;x) & \text{if} \quad n = 2m, \\ x \overline{T}_{n-1}(x) & \text{if} \quad n = 2m+1. \end{cases} \tag{73}$$

$$\overline{S}_n \begin{pmatrix} -1 & -2 \\ -1 & 1 \end{pmatrix} x = \begin{cases} \overline{S}_n(-1,1,-1,-2;x) & \text{if} \quad n = 2m, \\ x \overline{U}_{n-1}(x) & \text{if} \quad n = 2m+1. \end{cases} \tag{74}$$

However since we have

$$\int_{-1}^{1} W \begin{pmatrix} 1 & -2 \\ -1 & 1 \end{pmatrix} x \, dx = \int_{-1}^{1} \frac{dx}{x^2 \sqrt{1-x^2}} = +\infty \, , \tag{73.1}$$

$$\int_{-1}^{1} W \begin{pmatrix} -1 & -2 \\ -1 & 1 \end{pmatrix} x \, dx = \int_{-1}^{1} \frac{\sqrt{1-x^2}}{x^2} dx = +\infty \, , \tag{74.1}$$

they cannot fall into the half-trigonometric orthogonal polynomials category.

**Remark 1.** According to (50.4), the initial vector of the monic Chebyshev polynomials $\overline{U}_n(x)$ is $(p,q,r,s) = (-1,1,-3,0)$. If this vector is replaced in (18.1) a very simple case of three term recurrence relation (18) with $C_n(-1,1,-3,0) = -1/4$ is derived. A system of monic orthogonal polynomials that satisfies the relation

$$\overline{P}_{n+1}(x) = (x - \alpha_n) \overline{P}_n(x) - \beta_n \overline{P}_{n-1}(x) \, ; \quad \overline{P}_0(x) = 1 \, , \, n \in Z^+, \tag{75}$$

and has the property

$$\underset{n\to\infty}{Lim}\ \alpha_n = a \quad \text{and} \quad \underset{n\to\infty}{Lim}\ \beta_n = \frac{b^2}{4} > 0 \quad ; \quad a,b \in R, \tag{75.1}$$

is said to belong to the class $N(a,b)$. This class was introduced in detail by Nevai [18]. The monic polynomials corresponding to the conditions (75.1) are perturbations of $x \to b^n \overline{U}_n(\frac{x-a}{b})$. Now since

$$\underset{n\to\infty}{Lim}\ C_n\begin{pmatrix} -3 & 2 \\ -1 & 1 \end{pmatrix} = -\frac{1}{4} \quad \text{and} \quad \underset{n\to\infty}{Lim}\ C_n\begin{pmatrix} -5 & 2 \\ -1 & 1 \end{pmatrix} = -\frac{1}{4}, \tag{76}$$

the defined polynomials $\overline{X}_n(x)$ and $\overline{Y}_n(x)$ belong to the class $N(0,1)$.

## 4. Second subclass, Generalized Hermite polynomials (GHP)

The GHP were first introduced by Szego who gave a second order differential equation for these polynomials [22, problem 25] as almost the same form as we will give in this section. These polynomials can be characterized by using a direct relationship between them and Laguerre orthogonal polynomials [6]. Of course, there are some other approaches for this matter, see e.g. [8]. Because of this, it is better to only point to the main properties of GHP in terms of the obtained properties of $\overline{S}_n(p,q,r,s;x)$.

**Initial vector**

$$(p,q,r,s) = (0,1,-2,2a). \tag{77}$$

**Explicit form of monic GHP**

$$\overline{S}_n\begin{pmatrix} -2 & 2a \\ 0 & 1 \end{pmatrix}x = (-1)^{[\frac{n}{2}]}(a+1-\frac{(-1)^n}{2})_{[\frac{n}{2}]} \sum_{k=0}^{[n/2]} \binom{[n/2]}{k} \left(\prod_{i=0}^{[n/2]-(k+1)} \frac{-2}{2i+(-1)^{n+1}+2+2a}\right) x^{n-2k}. \tag{78}$$

**Recurrence relation of monic GHP**

$$\overline{S}_{n+1}(x) = x\overline{S}_n(x) + C_n\begin{pmatrix} -2 & 2a \\ 0 & 1 \end{pmatrix}\overline{S}_{n-1}(x) \quad ; \quad \overline{S}_0(x) = 1, \ \overline{S}_1(x) = x, \ n \in N, \tag{79}$$

where

$$C_n\begin{pmatrix} -2 & 2a \\ 0 & 1 \end{pmatrix} = -\frac{1}{2}n - \frac{1-(-1)^n}{2}a. \tag{79.1}$$

**Orthogonality relation**

$$\int_{-\infty}^{\infty} x^{2a} e^{-x^2} \overline{S}_n\begin{pmatrix} -2 & 2a \\ 0 & 1 \end{pmatrix}x\overline{S}_m\begin{pmatrix} -2 & 2a \\ 0 & 1 \end{pmatrix}x\,dx = \left(\frac{1}{2^n}\prod_{i=1}^{n}(1-(-1)^i)a+i\right)\Gamma(a+\frac{1}{2})\,\delta_{n,m}. \tag{80}$$

(80) shows that the orthogonality relation is valid for $a+1/2 > 0$ and $(-1)^{2a} = 1$.

**Differential equation:** $\Phi_n(x) = S_n(0,1,-2,2a;x)$.

$$x^2 \Phi_n''(x) - 2x(x^2 - a) \Phi_n'(x) + \left(2n x^2 + ((-1)^n - 1)a\right) \Phi_n(x) = 0 \ . \tag{81}$$

Finally it is worthy to point out that since the leading coefficient of Hermite polynomials $H_n(x)$ is $2^n$, the following equality holds

$$H_n(x) = 2^n \overline{S}_n\!\left(\begin{array}{cc|c} -2 & 0 \\ 0 & 1 \end{array} x\right) . \tag{82}$$

**5. Third subclass, A finite class of symmetric orthogonal polynomials with weight function $x^{-2a}(1+x^2)^{-b}$ on $(-\infty, \infty)$.**

According to Favard theorem, if the condition $-C_{n+1}(p,q,r,s) > 0$ holds only for a finite number of positive integers, i.e. for $n = 0,1,...,N$ then the related polynomials class would be finitely orthogonal. This note helps us obtain some new classes of finite symmetric orthogonal polynomials, which are special sub-cases of $\overline{S}_n(p,q,r,s;x)$ and can be indicated by it directly. To derive the first finite sub-class, let us compute the logarithmic derivative of the weight function $W(x) = x^{-2a}(1+x^2)^{-b}$ as

$$\frac{W'(x)}{W(x)} = \frac{-2(a+b)x^2 - 2a}{x(1+x^2)} \ . \tag{83}$$

If the above fraction is compared with the logarithmic derivative of the basic weight function $W(p,q,r,s;x)$ then we get

$$(p,q,r,s) = (1,1,-2a-2b+2,-2a) \ , \tag{84}$$

which is in fact the initial vector corresponding to the first finite sub-class of symmetric orthogonal polynomials. Hence, if (84) is replaced in (15) the explicit form of polynomials is derived as (85)

$$S_n\!\left(\begin{array}{cc|c} -2a-2b+2, & -2a \\ 1, & 1 \end{array} x\right) = \sum_{k=0}^{[n/2]} \binom{[\frac{n}{2}]}{k} \left(\prod_{i=0}^{[\frac{n}{2}]-(k+1)} \frac{2i + 2[n/2] + (-1)^{n+1} + 2 - 2a - 2b}{2i + (-1)^{n+1} + 2 - 2a}\right) x^{n-2k} \ .$$

Replacing (84) in the basic recurrence relation (18) also gives (86)

$$\overline{S}_{n+1}(x) = x \overline{S}_n(x) + C_n\!\left(\begin{array}{cc} -2a-2b+2, & -2a \\ 1, & 1 \end{array}\right) \overline{S}_{n-1}(x) \ ; \ \overline{S}_0(x) = 1, \ \overline{S}_1(x) = x, \ n \in N,$$

where

$$C_n\!\left(\begin{array}{cc} -2a-2b+2, & -2a \\ 1, & 1 \end{array}\right) = \frac{(n-(1-(-1)^n)a)(n-(1-(-1)^n)a-2b)}{(2n-2a-2b+1)(2n-2a-2b-1)} \ . \tag{87}$$

Therefore, the orthogonality relation

$$\int_{-\infty}^{\infty} \frac{x^{-2a}}{(1+x^2)^b} \bar{S}_n\!\left(\begin{array}{cc}-2a-2b+2, & -2a \\ 1, & 1\end{array}\bigg| x\right) \bar{S}_m\!\left(\begin{array}{cc}-2a-2b+2, & -2a \\ 1, & 1\end{array}\bigg| x\right) dx =$$
$$\left((-1)^n \prod_{i=1}^{n} C_i\!\left(\begin{array}{cc}-2a-2b+2, & -2a \\ 1, & 1\end{array}\right)\right) \frac{\Gamma(b+a-1/2)\Gamma(-a+1/2)}{\Gamma(b)} \delta_{n,m} \quad (88)$$

is valid if and only if $-C_{n+1}\!\left(\begin{array}{cc}-2a-2b+2 & -2a \\ 1 & 1\end{array}\right) > 0$ ; $\forall n \in Z^+$; $b+a > 1/2$; $a < 1/2$ and $b > 0$.

Here is a good position to explain how we can determine the parameters conditions to be established the orthogonality property (88). For this purpose, there is an interesting technique. Let us consider the differential equation of polynomials (85) using the initial vector (84) as

$$x^2(x^2+1)\Phi_n''(x) - 2x((a+b-1)x^2+a)\Phi_n'(x) + \left(n(2a+2b-(n+1))x^2 + (1-(-1)^n)a\right)\Phi_n(x) = 0. \quad (89)$$

If the latter equation is written in self-adjoint form, then according to theorem 1 the following term must vanish, i.e.

$$[x^{-2a}(1+x^2)^{-b+1}(\Phi_n'(x)\Phi_m(x) - \Phi_m'(x)\Phi_n(x))]_{-\infty}^{\infty} = 0 . \quad (90)$$

On the other hand, since $\Phi_n(x)$ is a polynomial of degree $n$, so

$$Max \deg(\Phi_n'(x)\Phi_m(x) - \Phi_m'(x)\Phi_n(x)) = n+m-1 . \quad (91)$$

Consequently from (90) and (91) we must have

$$-2a+2(-b+1)+n+m-1 \le 0 , \quad (92)$$

which gives the following result

$$\begin{cases}-2a-2b+n+m+1 \le 0, \\ a < 1/2, \; b > 0.\end{cases} \quad (93)$$

In other words, (88) holds if $m,n = 0,1,\ldots, N \le a+b-1/2$ where $N = Max\{m,n\}$ and $a < 1/2$; $b > 0$.

**Corollary 1**. The finite polynomial set $\{S_n(1,1,-2a-2b+2,-2a;x)\}_{n=0}^{n=N}$ is orthogonal with respect to the weight function $x^{-2a}(1+x^2)^{-b}$ on $(-\infty,\infty)$ if and only if $N \le a+b-1/2$, $a < 1/2$ and $b > 0$.

Note that the explained technique can similarly be applied for the first and second sub-classes of $S_n(p,q,r,s;x)$ i.e. GUP and GHP.

## 6. Fourth subclass, A finite class of symmetric orthogonal polynomials with weight function $x^{-2a}e^{-1/x^2}$ on $(-\infty, \infty)$.

Similar to the first finite sub-class, one can compute the logarithmic derivative of the given weight function to get respectively

**Initial vector**
$$(p,q,r,s) = (1,0,-2a+2,2). \tag{94}$$

**Explicit form of polynomials**

$$S_n\begin{pmatrix} -2a+2 & 2 \\ 1 & 0 \end{pmatrix} x = \sum_{k=0}^{[n/2]} \binom{[n/2]}{k} \left( \prod_{i=0}^{[\frac{n}{2}]-(k+1)} \frac{2i + 2[n/2] + (-1)^{n+1} + 2 - 2a}{2} \right) x^{n-2k}. \tag{95}$$

**Recurrence relation of monic polynomials**

$$\bar{S}_{n+1}(x) = x\bar{S}_n(x) + C_n\begin{pmatrix} -2a+2 & 2 \\ 1 & 0 \end{pmatrix} \bar{S}_{n-1}(x) \;;\quad \bar{S}_0(x) = 1,\; \bar{S}_1(x) = x,\; n \in N, \tag{96}$$

where

$$C_n\begin{pmatrix} -2a+2 & 2 \\ 1 & 0 \end{pmatrix} = \frac{-2(-1)^n(n-a) - 2a}{(2n-2a+1)(2n-2a-1)}. \tag{96.1}$$

**Orthogonality relation**

$$\int_{-\infty}^{\infty} x^{-2a} e^{-\frac{1}{x^2}} \bar{S}_n\begin{pmatrix} -2a+2 & 2 \\ 1 & 0 \end{pmatrix} x \bar{S}_m\begin{pmatrix} -2a+2 & 2 \\ 1 & 0 \end{pmatrix} x \, dx$$
$$= \left( (-1)^n \prod_{i=1}^{n} C_i\begin{pmatrix} -2a+2 & 2 \\ 1 & 0 \end{pmatrix} \right) \Gamma(a - \frac{1}{2}) \delta_{n,m}. \tag{97}$$

But (97) is valid if

$$[x^{2-2a} \exp(-\frac{1}{x^2})(\Phi'_n(x)\Phi_m(x) - \Phi'_m(x)\Phi_n(x))]_{-\infty}^{\infty} = 0, \tag{98}$$

or equivalently

$$2 - 2a + n + m - 1 \leq 0 \Leftrightarrow N \leq a - \frac{1}{2} \;;\quad N = Max\{m,n\}. \tag{99}$$

**Corollary 2.** The finite polynomial set $\{S_n(1,0,-2a+2,2;x)\}_{n=0}^{n=N}$ is orthogonal with respect to the weight function $x^{-2a}e^{-1/x^2}$ on $(-\infty,\infty)$ if and only if $N \leq a - 1/2$.

**Differential equation:** $\Phi_n(x) = S_n(1,0,-2a+2,2;x)$.

$$x^4 \Phi''_n(x) + 2x((1-a)x^2 + 1)\Phi'_n(x) - \left(n(n+1-2a)x^2 + 1 - (-1)^n\right)\Phi_n(x) = 0. \tag{100}$$

**7. A unified approach for the classification of BCSOP**

As we observed in the previous sections, each four introduced sub-classes of symmetric orthogonal polynomials were determined by $\bar{S}_n(p,q,r,s;x)$ directly and it was only sufficient to obtain the initial vector corresponding to them. On the other hand, it is clear that the orthogonality interval of the sub-classes, other than first one (GUP), are all infinite, i.e. $(-\infty,\infty)$. Hence, applying a linear transformation, say $x = wt + v$, preserves the orthogonality interval. For example, if $x = wt + v/w$ in (80), the orthogonality interval will not change and consequently a more extensive class with weight function $(w^2 t + v)^{2a} \exp(-w^2 t^2 - 2vt)$ will be derived on the interval $(-\infty,\infty)$. However, it is important to know that the latter weight corresponds to the class of orthogonal polynomials $S_n(p,q,r,s;wx+v/w)$. Therefore, only by having the initial vector we can have access to all other standard properties and design a unified approach for the cases that may occur. In other words, if one can obtain the parameters $(p,q,r,s)$ by referring to the initial data, such as a given three term recurrence relation, a given weight function and so on, then all other properties will be derived straightforwardly. Here, we consider two special sub-cases of this approach.

**7.1. How to find the parameters $p,q,r,s$ if a special case of the main weight function $W(p,q,r,s;x)$ is given?**

By referring to third and fourth orthogonal sub-classes, it is easy to find out that the best way for deriving $p,q,r,s$ is to compute the logarithmic derivative $W'(x)/W(x)$ and then equate the pattern with (21.1). The following examples will clarify this matter.

**Example 1.** The weight functions

i)      $W_1(x) = -x^6 + 4x^4$                     ;    $-2 \leq x \leq 2$

ii)     $W_2(x) = (16x^2 - 8x + 1)\exp(2x(1-2x))$ ;    $-\infty < x < \infty$

iii)    $W_3(x) = (2x+1)^2 (2x^2 + 2x + 1)^{-5}$         ;    $-\infty < x < \infty$

and their orthogonality intervals are given. Find other standard properties such as explicit form of polynomials, orthogonality relation and ….

To solve the problem, it is only sufficient to find the initial vector corresponding to each given weight functions. For this purpose, if the logarithmic derivative of the first weight is computed, then we have

i)    $\dfrac{W_1'(x)}{W_1(x)} = \dfrac{6x^2 - 16}{x(x^2 - 4)} = \dfrac{(r - 2p)x^2 + s}{x(px^2 + q)} \Rightarrow (p,q,r,s) = (1, -4, 8, -16)$.

Hence the related monic polynomials are $\{\bar{S}_n(1,-4,8,-16;x)\}_{n=0}^{\infty}$. Note that these polynomials are orthogonal with respect to $x^4(4 - x^2)$ on $[-2, 2]$ for every $n$ and it is *not* necessary to know that they are the same as shifted GUP on $[-2,2]$, because they can explicitly and independently be expressed by $\bar{S}_n(p,q,r,s;x)$.

But, for the weight function $W_2(x)$ we have

ii) $W_2(x) = 4e^{\frac{1}{4}}(2x-\frac{1}{2})^2 \exp(-(2x-\frac{1}{2})^2) \Rightarrow \frac{1}{4}e^{-\frac{1}{4}}W_2(\frac{2t+1}{4}) = t^2 e^{-t^2} = W_2^*(t)$

$\rightarrow \frac{W_2^{*'}(t)}{W_2^*(t)} = \frac{-2t^2+2}{t} \Rightarrow (p,q,r,s) = (0,1,-2,2)$

Hence the related orthogonal polynomials are as $\left\{S_n(0,1,-2,2; 2x-\frac{1}{2})\right\}_{n=0}^{\infty}$.

iii) $W_3(x) = 2^5 \frac{(2x+1)^2}{(1+(2x+1)^2)^5} \Rightarrow 2^{-5}W_3(\frac{t-1}{2}) = \frac{t^2}{(1+t^2)^5} = W_3^*(t)$

$\rightarrow \frac{W_3^{*'}(t)}{W_3^*(t)} = \frac{-8t^2+2}{t(t^2+1)} \Rightarrow (p,q,r,s) = (1,1,-6,2)$

Therefore, by noting the orthogonality relation of the third sub-class of BCSOP, the finite set $\{S_n(1,1,-6,2; 2x+1)\}_{n=0}^{n=3}$ is orthogonal with respect to $W_3(x)$ on $(-\infty,\infty)$ and the upper bound of this set has been determined based on the condition $N \leq a+b-1/2$ for $b=5$ and $a=-1$.

## 7.2. How to find the parameters $p,q,r,s$ if a special case of the main three-term recurrence equation (18) is given?

In general, there are two ways to determine the special case of $\bar{S}_n(p,q,r,s;x)$ corresponding to a given three-term recurrence equation. The first way is to directly compare the given recurrence equation with (18). This leads to a system of polynomial equations in terms of the four parameters $p,q,r,s$. The second way is to equate the first four terms of each two recurrence equations together, which leads to a polynomial system with 4 equations and 4 unknowns $p,q,r$ and $s$ respectively. The following example will better illustrate these methods.

**Example 2.** If the recurrence equation

$$\bar{S}_{n+1}(x) = x\bar{S}_n(x) - 2\frac{6+(-1)^n(n-6)}{(2n-11)(2n-13)}\bar{S}_{n-1}(x) \quad ; \quad \bar{S}_0(x)=1 \text{ and } \bar{S}_1(x)=x,$$

is given, then find its explicit polynomial solution, differential equation of polynomials, the related weight function and finally orthogonality relation of polynomials.

**Solution.** If the above recurrence equation is directly compared with the main equation (18) and subsequently (18.1), then one can obtain the values $(p,q,r,s) = (1,0,-10,2)$. Hence, the explicit solution of above recurrence equation is the polynomials $\bar{S}_n(1,0,-10,2;x)$ and therefore consequently their differential equation is found as

$$x^4 \Phi_n''(x) + x(-10x^2 + 2)\Phi_n'(x) - \left(n(n-11)x^2 + 1 - (-1)^n\right)\Phi_n(x) = 0.$$

Moreover, by replacing the initial vector in the main weight function (21.1) as

$$W\begin{pmatrix} -10 & 2 \\ 1 & 0 \end{pmatrix} x = \exp\left(\int \frac{-12x^2 + 2}{x^3} dx\right) = x^{-12} e^{-\frac{1}{x^2}},$$

one can find out that the related polynomials are a particular case of the fourth introduced sub-class. Consequently we have

$$\int_{-\infty}^{\infty} x^{-12} e^{-\frac{1}{x^2}} \bar{S}_n\begin{pmatrix} -10 & 2 \\ 1 & 0 \end{pmatrix} x \bar{S}_m\begin{pmatrix} -10 & 2 \\ 1 & 0 \end{pmatrix} x dx = \left((-1)^n \prod_{i=1}^{i=n} C_i\begin{pmatrix} -10 & 2 \\ 1 & 0 \end{pmatrix}\right) \Gamma(\tfrac{11}{2}) \delta_{n,m} \Leftrightarrow m, n \leq 5$$

**Second method.** If the given recurrence relation is only expanded for $n = 2, 3, 4, 5$ and then equated with (17.1), the following system will be derived

$$\bar{S}_2(x) = x^2 - \frac{2}{9} = x^2 + \frac{q+s}{p+r} \Rightarrow \frac{q+s}{p+r} = -\frac{2}{9} \;;\; \bar{S}_3(x) = x\,\bar{S}_2(x) - \frac{4}{63}\bar{S}_1(x) = x^3 + \frac{3q+s}{3p+r}x$$

$$\bar{S}_4(x) = x\,\bar{S}_3(x) - \frac{18}{35}\bar{S}_2(x) = x^4 + 2\frac{3q+s}{5p+r}x^2 + \frac{(3q+s)(q+s)}{(5p+r)(3p+r)} \quad and \quad 2\frac{5q+s}{7p+r} = -\frac{4}{3}.$$

Solving this system again results that $(p, q, r, s) = (1, 0, -10, 2)$.

**8. Conclusion.** By using the extended Sturm-Liouville theorem for symmetric functions, one can define a generic second order differential equation having a basic polynomial solution with four free parameters. This solution satisfies a generic orthogonality relation whose weight function corresponds to an analogue of Pearson distributions. In other words, there are four special cases of the dual symmetric distributions family that can respectively be considered as the weight functions of four introduced sub-classes of BCSOP. In this way, the following table shows the explicit forms of the mentioned sub-classes in terms of $S_n(p, q, r, s; x)$ as well as their weight functions, kind of polynomials (finite or infinite), orthogonality interval and finally constraint on the parameters.

**Table 1: Four special sub-cases of $S_n(p, q, r, s; x)$**

| Definition | Weight function | Interval & Kind | Parameters Constraint |
|---|---|---|---|
| $S_n\begin{pmatrix} -2a-2b-2, & 2a \\ -1, & 1 \end{pmatrix} x$ | $x^{2a}(1-x^2)^b$ | $[-1,1]$, Infinite | $a > -1/2$ $b > -1$ |
| $S_n\begin{pmatrix} -2, & 2a \\ 0, & 1 \end{pmatrix} x$ | $x^{2a} e^{-x^2}$ | $(-\infty, \infty)$, Infinite | $a > -\frac{1}{2}$ |
| $S_n\begin{pmatrix} -2a-2b+2, & -2a \\ 1, & 1 \end{pmatrix} x$ | $\dfrac{x^{-2a}}{(1+x^2)^b}$ | $(-\infty, \infty)$, Finite | $\begin{cases} N \le a+b-1/2 \\ a < 1/2,\ b > 0 \end{cases}$ |
| $S_n\begin{pmatrix} -2a+2, & 2 \\ 1, & 0 \end{pmatrix} x$ | $x^{-2a} e^{-\frac{1}{x^2}}$ | $(-\infty, \infty)$, Finite | $N \le a - \dfrac{1}{2}$ |

Note that since all weights in above table are even functions, the condition $(-1)^{2a} = 1$ must always be satisfied by noting the constraint of parameters for each introduced weight functions. Therefore, we can also consider them in the forms $|x|^{2a}(1-x^2)^b$, $|x|^{2a} e^{-x^2}$, $|x|^{-2a}(1+x^2)^{-b}$ and $|x|^{-2a} e^{-1/x^2}$ respectively.

**Remark 2.** The differential equation (10) can be extended to a generic operator equation of the form

$$x^2(p^*x^2 + q^*)L^2 P_n(x) + x(r^*x^2 + s^*)L P_n(x) - \left(n(r^* + (n-1)p^*)x^2 + (1-(-1)^n)s^*/2\right) P_n(x) = 0$$

where the linear operator $L$, known as Hahn's operator [6, p.159], is defined by

$$L(f(x)) = \frac{f(qx+w) - f(x)}{(q-1)x + w} \ ;\ q, w \in \mathbf{R}.$$

Clearly the ordinary derivative operator, which is the limiting case $w = 0$ and $q \to 1$, corresponds to equation (10). Hence, two other cases can be considered. First case is the difference operator $\Delta$ for $q = w = 1$ (or equivalently $\nabla$ for $q = -w = 1$). Similar to the continuous case, since the extension of Sturm-liouville problems is possible for discrete variables there is a basic class of discrete orthogonal polynomials that satisfies a generic difference equation and possesses four special sub-classes [17]. The second case, i.e. q-difference operator $D_q$ for $w = 0$ can be studied and investigated individually.

**Remark 3.** The extended Sturm-liouville theorem for symmetric functions may essentially open a new perspective. For instance, the following differential equation

$$x^2(ax^\lambda + b) y'' + x(cx^\lambda + d) y' + \left(\alpha_n^* x^\lambda + \beta_n^*\right) y = 0,$$

is a generalization of differential equation (10) for $\lambda = 2$. According to [16] this equation has a general solution, which is orthogonal with respect to an even weight function for specific sequences of $\alpha_n^*$ and $\beta_n^*$.

**References.**

[1] Abramowitz, M. and Stegun, I. A. (1972), Handbook of mathematical functions with formulas, graphs, and mathematical tables, 9[th] printing, New York: Dover

[2] Al-Salam, W., Allaway, W.R. and Askey, R. (1984), Sieved ultraspherical polynomials, Trans. Amer. Math. Soc., **284**, pp. 39–55

[3] Arfken, G. (1985), Mathematical methods for physicists, Academic Press Inc

[4] Askey, R. (1984), Orthogonal polynomials old and new, and some combinatorial connections,Enumeration and Design (D.M. Jacson and S.A. Vanstone, eds.), Academic Press, New York, pp. 67-84

[5] Charris, J.A. and Ismail, M.E.H. (1986), On sieved orthogonal polynomials II: Random walk polynomials, Canad. J. Math. **38**, pp. 397 – 415

[6] Chihara, T.S. (1978), Introduction to Orthogonal Polynomials, Gordon & Breach, New York

[7] Cover, T.M. and Thomas, J.A. (1991), Elements of information theory, Wiley-Interscience, New York

[8] Dette, H. (1996), Characterizations of generalized Hermite and sieved ultraspherical polynomials, Proc. Amer. Math. Soc, **384**, pp. 691 – 711

[9] Fisher, R.A. (1925), Theory of statistical estimation, Proc. Cambridge Philos. Soc. **22** pp. 700 -725

[10] Fox, L. and Parker, I.B. (1968), Chebyshev polynomials in numerical Analysis, London, England, Oxford University Press

[11] Ismail, M.E.H. (1986), On sieved orthogonal polynomials III: Orthogonality on several intervals, Trans. Amer. Math. Soc. **294**, pp. 89 – 111

[12] Koepf, W. (1988), Hypergeometric Summation, Braunschweig/Wiesbaden: Vieweg

[13] Konoplev, V.P. (1965), The asymptotic behavior of orthogonal polynomials at one-sided singular points of weighting functions (algebraic singularities), Soviet Math. Doklady. **6**, pp. 223-227

[14] Masjed-Jamei, M. (2002), Three finite classes of hypergeometric orthogonal polynomials and their application in functions approximation, Integral Transforms and Special Functions, **13**, pp. 169-190

[15] Masjed-Jamei, M. (2004), Classical orthogonal polynomials with weight function $((ax+b)^2 + (cx+d)^2)^{-p} \exp(q \arctan \frac{ax+b}{cx+d})$ ; $x \in (-\infty, \infty)$ and a generalization of T and F distributions, Integral Transforms and Special Functions, **15**, pp. 137-153

[16] Masjed-Jamei, M., (2007), A generalization of classical symmetric orthogonal functions using a symmetric generalization of Sturm-Liouville problems, Integral Transforms and Special Functions, **18**, pp. 871-883.

[17] Masjed-Jamei, M., Area, I., (2013), A basic class of symmetric orthogonal polynomials of a discrete variable, Journal of Mathematical Analysis and Applications, **399**, pp. 291-305.

[18] Nevai, p. (1979), Orthogonal polynomials, Memoirs Amer. Math. Soc., Vol. 213, Amer. Mat. Soc., providence, R.I.

[19] Nikiforov, A. F. and Uvarov, V. B. (1988), Special Functions of Mathematical Physics, Basel-Boston: Birkhäuser

[20] Rivlin,T. (1990), Chebyshev polynomials: From approximation theory to algebra and number theory, 2th edition, John Wiley and Sons, New York

[21] Rogers, L.J. (1895), Third memoir on the expansion of certain infinite products, Proc. London Math. Soc. **26**, pp. 15 – 32

[22] Szegö, G. (1975), Orthogonal polynomials. American Mathematical Society Colloquium Publications, Vol. 23, Providence, RI